\documentclass[11pt]{article}
\usepackage{epsf,amssymb,latexsym,amsmath}
\usepackage{graphicx}
\usepackage{tikz}
\usepackage[autostyle=true]{csquotes}

\newcommand{\proof}{\noindent{\bf Proof.\ }}
\newcommand{\qed}{\hfill $\square$ \bigskip}

\newtheorem{theorem}{\bf Theorem}[section]
\newtheorem{corollary}[theorem]{\bf Corollary}

\newtheorem{conjecture}[theorem]{\bf Conjecture}

\begin{document}

\title{A note on Vizing's conjecture}

\author{
Simon \v Spacapan\footnote{ University of Maribor, FME, Smetanova 17,
2000 Maribor, Slovenia and IMFM, Jadranska 19, 1000 Ljubljana, Slovenia. e-mail: simon.spacapan@um.si.
}}
\date{\today}

\maketitle

\begin{abstract}
Let $\gamma(G)$ denote the domination number of graph $G$.   Let $G$ and $H$ be   graphs and $G\Box H$ their Cartesian product. 
 For  $h\in V(H)$   define   
$G_h=\{(g,h)\,|\,g\in V(G)\}$ and call this set a $G$-layer of $G\Box H$. 
We prove the following special case of Vizing's conjecture. 
 Let $D$ be a dominating set of $G\Box H$. 
If there exist minimum dominating sets  $D_1$ and $D_2$ of $G$ such that 
for every $h\in V(H)$, the projection of $D\cap G_h$ to $G$ is contained in $D_1$ or $D_2$, then $|D|\geq \gamma(G)\gamma(H)$.
\end{abstract}





\section{Introduction}
In 1968 Vizing conjectured that for any graphs $G$ and $H$ we have $\gamma(G\Box H)\geq \gamma(G)\gamma(H)$, see \cite{vizing}.  
The motivation for  this article is derived from a well known result, obtained by Clark and Suen in \cite{clarksuen}, which asserts 
that $\gamma(G\Box H)\geq \frac 12 \gamma(G)\gamma(H)$. 
Note however that 
this result was improved by S.~Zerbib in \cite{zerbib} (see also \cite{brezar}) where it is proved that 
$\gamma(G\Box H)\geq \frac 12 \gamma(G)\gamma(H) +\frac 12 \max\{ \gamma(G),\gamma(H)\}$. The proof of the result given in  \cite{clarksuen} is outlined and discussed below (we use the notation from  \cite{imklarall}).

Let $G$ and $H$ be graphs and   $D$ a dominating set of $G\Box H$. 
We say that a vertex $(u,v)\in V(G\Box H)$ is {\em vertically dominated} by $D$ if there exists a vertex $v'\in N_H[v]$ such that $(u,v')\in D$. 

 Let $\{g_1,\ldots,g_k\}$ be a minimum dominating set of $G$, and let
 $\{\pi_1,\ldots, \pi_k\}$ be any partiton of $V(G)$ such that $\pi_i\subseteq N[g_i]$, where $1\leq i\leq k=\gamma(G)$. For   $i\in [k]$    and   $h\in V(H)$
we call the set $\pi_i\times\{h\}$ {\em a cell}. A cell is {\em vertically undominated} if no vertex from this cell is vertically dominated by $D$. Define $H_i=\pi_i\times V(H)$.

For $i\in [k]$, let $n_i$ be the number of vertically undominated cells contained in $H_i$. 
For $h\in V(H)$, let $m_h$ be the number of vertically undominated cells contained in $G_h$. The result of 
Clark and Suen follows from combining \eqref{prva} and \eqref{druga} given below, where \eqref{prva} follows from 
$|D\cap H_i|+n_i\geq \gamma(H)$ for every $i\in [k]$, and \eqref{druga} follows from $|D\cap G_h|\geq m_h$ for every $h\in V(H)$
(note that  $|D\cap G_h|< m_h$ would contradict our assumption that $\gamma(G)=k$).
\begin{equation}\label{prva}
|D|+\sum_{i=1}^k n_i=\sum_{i=1}^k (|D\cap H_i|+n_i)\geq \gamma(G)\gamma(H)
\end{equation}

\begin{equation}\label{druga}
|D|=\sum_{h\in V(H)}|D\cap G_h|\geq \sum_{h\in V(H)}m_h=\sum_{i=1}^k n_i
\end{equation}

Thus we have  $\gamma(G\Box H)\geq \frac 12 \gamma(G)\gamma(H)$,  moreover we note that $\gamma(G\Box H)= \frac 12 \gamma(G)\gamma(H)$ is 
possible only if there is an equality in   \eqref{prva} and \eqref{druga}. 
In particular this means that $|D\cap G_h|=m_h$ for every $h\in V(H)$. 
In this article we ask the question: \enquote{what can we say about (the size of) $D$ if we assume that $|D\cap G_h|=m_h$ for every $h\in V(H)$?} 

Now suppose that $|D\cap G_h|=m_h$ and that cells $\pi_i\times \{h\}$ for $i\in I$ are vertically undominated cells, where $|I|=m_h$. 
Then $p_G(D\cap G_h)\cup \{g_i\,|\,i\notin I\}$ is a minimum dominating set of $G$, and consequently  $p_G(D\cap G_h)$ is contained in a minimum dominating 
set of $G$. 
Hence, we can formulate the above question slightly more general as follows: \enquote{what can we say about (the size of) $D$, if we assume that for every 
$h\in V(H)$ the set $ p_G(D\cap G_h)$ is contained in a minimum dominating set of $G$ ?}  In view of this question the following weaker version of Vizing's conjecture is the object of our interest.

\begin{conjecture}
Let $G$ and $H$ be graphs and let $D$ be a dominating set of $G\Box H$. If for every $h\in V(H)$ there is a minimum dominating set 
$D_h$ of $G$ such that $p_G(D\cap G_h)\subseteq D_h$, then $|D|\geq \gamma(G)\gamma(H)$.
\end{conjecture}

In this paper we prove (see Corollary~\ref{cor}) a special case of the above conjecture.  However, in order to pose our results more generally, 
we formulate our main result, Theorem~\ref{main},  in   terms of minimal dominating sets, which are defined in the sequel. 

Let $G$ be a graph and $X\subseteq V(G)$. A {\em minimal dominating set of $G$ containing $X$} is a dominating  set $D$ of $G$, such that $D'$ is not a dominating set of 
$G$ whenever $X\subseteq D'\subset D$. 
Let $M_G(X)$ be the set of all minimal dominating sets of $G$ containing $X$. Note that $M_G(X)$ is nonempty for every $X\subseteq V(G)$. 

Let $S\subseteq V(G)$ be a set and suppose that $y\notin S$ is a neighbor of $x\in S$. We say that $y$ is a {\em private neighbor} of $x$ with respect to $S$
if $x$ is the only vertex of $S$ which is adjacent to $y$, that is $N(y)\cap S=\{x\}$. We will say that \enquote{$y$ is a private neighbor of $x\in S$}, and this shall be understood that 
$y$ is a private neighbor of $x\in S$ with respect to $S$ 
(unless otherwise stated, private neighbor of a vertex 
is always with respect to the set from which this vertex is taken). 

Let $G$ and $H$ be graphs and $G\Box H$ their Cartesian product.
For an $h\in V(H)$ the set $G_h=\{(g,h)\,|\,g\in V(G)\}$ is called a {\em G-layer}, and for any $g\in V(G)$  the set 
$H_g=\{(g,h)\,|\,h\in V(H)\}$ is called an {\em H-layer}. The mapping $p_G:V(G\Box H)\rightarrow V(G)$ defined by 
$p_G:(g,h)\mapsto g$ is called the {\em projection} from $V(G\Box H)$ to $V(G)$. Let $D$ be a dominating set of $G\Box H$. 
We say that a vertex $(u,v)\in V(G\Box H)$ is {\em horizonatally dominated} by $D$ if there exists a vertex $u'\in N_G[u]$ such that $(u',v)\in D$. 
For any positive integer $n$ we set $[n]=\{1,\ldots,n\}$. The complement of set $X$ is denoted by $\bar X$.

\section{Results}

Let $G$ and $H$ be graphs and let $D$ be a dominating set of $G\Box H$.
 In the following theorem we assume that there exists a fixed dominating set $S$ of $G$ such that for every   $h\in V(H)$ the set $S$ is a minimal dominating set containing  $p_G(D\cap G_h)$.

\begin{theorem} \label{1}
Let $G$ and $H$ be graphs and let $D$ be a dominating set of $G\Box H$. If there is a dominating set $S$ of $G$ such that  for every $h\in V(H)$  we have $S\in M_G(p_G(D\cap G_h))$, then $|D|\geq \gamma(G)\gamma(H)$.
\end{theorem}

\proof
Let $S\subseteq V(G)$ be   such that for every $h\in V(H)$ the set $S$ is a minimal dominating set containing $p_G(D\cap G_h)$. 
Let $S_1$ be the set of vertices $x\in S$ that have a private neighbor (with respect to $S$), and define $S_2=S\setminus S_1$.  
If $y\notin S$ is a private neighbor of $x\in S$, then 
$(y,h)$ is dominated by $(x,h)$, because $H_y\cap D=\emptyset$. It follows that $S_1\times V(H)\subseteq D$.  If $x\in S_2$ and $(x,h)\notin D$ for some 
$h\in V(H)$, then $x$ is not adjacent to any vertex of $S$, for otherwise $S$ is not a minimal dominating set containing $p_G(D\cap G_h)$ 
(note that $S-x$ dominates $G$ and $ p_G(D\cap G_h)\subseteq S-x$). It follows that 
every vertex in $H_x,x\in S_2$ is vertically dominated, equivalently $|D\cap H_x|\geq \gamma (H)$ for every $x\in S_2$. Recall that $S_1\times V(H)\subseteq D$ and so 
$|D\cap H_x|\geq \gamma (H)$ for every $x\in S$. Since $S$ is a dominating set of $G$ we have $|D|\geq \gamma (G)|S|\geq  \gamma(G)\gamma(H)$.
\qed

 In the following theorem we assume that there exist two dominating sets $S_1$ and $S_2$ of $G$ such that for every   $h\in V(H)$, $S_1$ or $S_2$ is a minimal dominating set containing  $p_G(D\cap G_h)$.

\begin{theorem} \label{main}
Let $G$ and $H$ be graphs and let $D$ be a dominating set of $G\Box H$. If there exist dominating sets $S_1$ and $S_2$ of $G$ such that  for every $h\in V(H)$  we have $S_i\in M_G(p_G(D\cap G_h))$ for some  $i\in [2]$, then $|D|\geq \gamma(G)\gamma(H)$.
\end{theorem}

\proof
Assume that $S_1$ and $S_2$ are dominating sets of $G$ such that for every $h\in V(H),$  $S_1$ or $S_2$ is a minimal dominating set containing 
 $p_G(D\cap G_h)$. We may also assume that $|S_1|\leq |S_2|$. Let $T_1$ be the set of $h\in V(H)$ such that $S_1$ is a minimal dominating set containing 
 $p_G(D\cap G_h)$, and define $T_2=V(H)\setminus T_1$. By Theorem \ref{1} we may assume that $T_1$ and $T_2$ are nonempty.  
We state several properties of $D$.
\begin{itemize}
\item[(i)] 
If a vertex $x\in S_i$, where $i\in [2]$, has a private neighbor (with respect to $S_i$) in $V(G)\setminus (S_1\cup S_2)$, then $\{x\}\times T_i\subseteq D$. 
\item[(ii)] If a vertex $x\in S_i$, where $i\in [2]$, has no private neighbors (with respect to $S_i$), then $\{x\}\times T_i\subseteq D$ or $x$ is not adjacent to any vertex in $S_i$. 
\end{itemize}
To prove (i) observe that no vertex in $H_y$, where  $y\in V(G)\setminus (S_1\cup S_2)$, is vertically dominated. 
So if $x\in S_i$ has a private neighbor $y\in V(G)\setminus (S_1\cup S_2)$, then for every $h\in T_i$, vertex $(y,h)$ is dominated by $(x,h)$, and so 
$(x,h)\in D.$ This proves (i).  

To prove (ii) suppose   that   $x\in S_i$ has no private neighbors and $x$ is adjacent to a vertex in $S_i$ for some $i\in [2]$. Then $S_i-x$ is a 
dominating set of $G$. Since $S_i$ is a minimal dominating set containing $p_G(D\cap G_h)$, we find that 
 $x\in p_G(D\cap G_h)$ for every $h\in T_i$. Hence, $\{x\}\times T_i\subseteq D$, which completes the proof of (ii). Note that (ii) 
implies that if $x\in S_i$ has no private neighbors, then every vertex in $\{x\}\times T_i$ is vertically dominated.

Let  $S_i=A_i\cup B_i\cup C_i$ be a partitioning of $S_i$, where $A_i,B_i,C_i$ for $i\in [2]$ are defined below.
\begin{itemize}
\item[(a)] $A_i$ is the set of vertices $x\in S_i$ that have a private neighbor in $V(G)\setminus (S_1\cup S_2)$.
\item[(b)]  $B_i$ is the set of vertices $x\in S_i$ that have no private neighbors in  $V(G)\setminus (S_1\cup S_2)$, but have a private neighbor in $S_j\setminus S_i$, where  $j\in [2]$ and $ j\neq i$.
\item[(c)]  $C_i$ is the set of vertices $x\in S_i$ that have no private neighbors. 
\end{itemize}

We claim that for every $x\in (S_1\cap S_2)\setminus (B_1\cap B_2)$, every vertex in $H_x$ is vertically dominated and so $|D\cap H_x|\geq \gamma (H)$. 
Note that $$(S_1\cap S_2)\setminus (B_1\cap B_2)= (A_1\cup B_1\cup C_1)\cap(A_2\cup B_2 \cup C_2)\setminus (B_1\cap B_2)\,.$$

Let us first consider $x\in A_1\cap A_2$. By (i) we get $\{x\}\times T_1\subseteq D$ and $\{x\}\times T_2\subseteq D$. Therefore 
$H_x\subseteq D$ and so every vertex in $H_x$ is vertically dominated. 

If $x\in A_1\cap B_2$, then $\{x\}\times T_1\subseteq D$ and 
$x$ has a private neighbor $y\in S_1\setminus S_2$. Each vertex $(y,h)$, where $h\in T_2$ is dominated by some vertex in 
$\{y\}\times T_1$ or by $(x,h)$. In other words, if $(x,h)\notin D$ then $(y,h)$ has a neighbor in $\{y\}\times T_1$ and so  
$(x,h)$ has a neighbor in $\{x\}\times T_1$. It follows that every vertex in 
$H_x$ is vertically dominated.  

If $x\in A_1\cap C_2$ then by (i) we have $\{x\}\times T_1\subseteq D$ and by (ii) every vertex in $\{x\}\times T_2$ is vertically dominated. Hence, all vertices in 
$H_x$ are vertically dominated. 

If $x\in B_1\cap A_2$ the arguments are analogous to the case when $x\in  A_1\cap B_2.$

If $x\in B_1\cap C_2$ then $x$ has a (private) neighbor $y$ in $S_2\setminus S_1$ and therefore, by (ii) and the fact that $x\in C_2$,  $\{x\}\times T_2\subseteq D$. Now for every $h\in T_1$, if $(x,h)$ is not adjacent to a vertex  in $\{x\}\times T_2$ then $(x,h)\in D$ 
(this follows from the fact that $(y,h)$ is not vertically dominated and $(x,h)$ is the only vertex which can horizontally dominate $(y,h)$). Hence 
every vertex in $H_x$ is vertically dominated. 

If $x\in C_1\cap A_2$ or $x\in C_1\cap B_2$ the arguments are analogous to cases $x\in A_1\cap C_2$ and $x\in B_1\cap C_2$, respectively. 

If $x\in C_1\cap C_2$ then, by (ii), every vertex in   $\{x\}\times T_i$ is vertically dominated, for $i\in [2]$. This proves the claim. 

Next we claim that for every $x\in B_1\cap B_2$ there exist vertices $y_x\in S_1\setminus S_2$ and $z_x\in S_2\setminus S_1$ such that 
$|(H_x\cup H_{y_x}\cup H_{z_x})\cap D|\geq 2\gamma (H)$ and so that both mappings 
$$x\mapsto y_x~~{\rm and}~~x\mapsto z_x$$
are one to one (injective). 

To prove the claim, for every $x\in B_1\cap B_2$ let $y_x\in S_1\setminus S_2$ be a private neighbor of $x$ (with respect to $S_2$). Note that the mapping  
$x\mapsto y_x$ is one to one because each vertex in $\bar S_2$ is a private neighbor of at most one vertex in $S_2$. 
Now observe that $y_x\notin B_1$ because each vertex in $B_1$ has a neighbor  $u$ in $S_2\setminus S_1$, and so $y_x\in B_1$ implies that  
$y_x$ is adjacent to $x\in S_2$ and $u\in S_2$, and therefore $y_x$ is not a private neighbor of $x$ (which contradicts our choice of $y_x$).  It follows that 
$y_x\in A_1$ or $y_x\in C_1$ and so each vertex in $\{y_x\}\times T_1$ is vertically dominated. Moreover, if a vertex $(y_x,h),h\in T_2$ is not 
vertically dominated, then it is horizontally dominated by $(x,h)$. It follows that $|((\{x\}\times T_2)\cup (\{y_x\}\times T_1))\cap D|\geq \gamma (H)$. 

For every $x\in B_1\cap B_2$ let $z_x\in S_2\setminus S_1$ be a private neighbor of $x$  (with respect to $S_1$).
We prove, similarly as above, that  $|((\{x\}\times T_1)\cup (\{z_x\}\times T_2))\cap D|\geq \gamma(H)$ which completes  the proof of the claim. 

Define $Q=\{y_x,z_x\,|\,x\in B_1\cap B_2\}$, and observe that $Q$ is disjoint with $B_1\cup B_2$ (we proved above that $y_x\notin B_1$, and we can similarly prove that $z_x\notin B_2$). We claim that there is a matching $M\subseteq E(G)$ such that 
\begin{itemize}
\item[(i)] for every edge $xy\in M$, $|(H_x\cup H_y)\cap D|\geq \gamma (H)$;
\item[(ii)] every edge of $M$ has at lest one endvertex in $(B_1\setminus S_2)\cup (B_2\setminus S_1)$; 
\item[(iii)] every vertex in  $(B_1\setminus S_2)\cup (B_2\setminus S_1)$ is saturated by $M$, and 
\item[(iv)] no edge of $M$ is incident to a vertex in $Q$.
\end{itemize}

To prove the claim let us first make the following observation.

\medskip
\noindent {\em Observation 0: If $b\in B_2\setminus S_1$ is a private neighbor of $a\in B_1\setminus S_2$ (with respect to $S_1$), then 
$a$ is  a private neighbor of $b$ (with respect to $S_2$). Moreover $b$ is the only neighbor of $a$ in $S_2\setminus S_1$, and 
$a$ is the only neighbor of $b$ in $S_1\setminus S_2$.}

\medskip
The observation follows from the fact that $b\in B_2$  has a private neighbor $c\in S_1\setminus S_2$ 
(with respect to $S_2$), and if $c\neq a$ then (by the definition) $b$ is not a private neighbor of $a$.  

 For every $x\in B_1\setminus S_2$  let $u_x\in S_2\setminus S_1$ be a private neighbor of $x$ (with respect to $S_1$), and observe that 
$u_x\notin Q$ (because every vertex in $Q$ is a private neighbor of a vertex in $B_1\cap B_2$, and each vertex in $S_2\setminus S_1$ is a private neighbor of at most one vertex in $S_1$). Moreover, for $x\in B_1\setminus S_2$ we also have $x\notin Q$ (recall that $Q$ is disjoint with $B_1\cup B_2$), and hence no edge  $xu_x$ 
has an endvertex in $Q$. 

If $u_x\in A_2\cup C_2$, then every vertex in $\{u_x\}\times T_2$ is vertically dominated. Moreover, since $u_x$ is a private neighbor of $x$, every vertex $(u_x,h),h\in T_1$ which is not vertically dominated is horizontally dominated by $(x,h)$. It 
follows that  $|(H_x\cup H_{u_x})\cap D|\geq \gamma (H)$. 

If $u_x\in B_2$, then by Observation 0, $x$ is a private neighbor of $u_x$. 
Hence, every vertex  $(x,h),h\in T_2$ which is not vertically dominated is horizontally dominated by $(u_x,h)$, and (as before)
every vertex $(u_x,h),h\in T_1$ which is not vertically dominated is horizontally dominated by $(x,h)$, 
so we conclude that also in this case   we have 
 $|(H_x\cup H_{u_x})\cap D|\geq \gamma (H)$. 

The reader may note that (the above discusssion proves that) the set of edges $xu_x, x\in B_1\setminus S_2$ is a matching that fulfills (i),(ii) and (iv), and saturates 
 $B_1\setminus S_2$. 

Let $R$ be the set of  $x'\in B_2\setminus S_1$, such that $x'\neq u_x$ for all  $x\in B_1\setminus S_2$, and 
note that (by Observation 0) no vertex in $R$ has a private neighbor in $B_1\setminus S_2$. 

 For every $x'\in R$,  let $w_{x'}\in S_1\setminus S_2$ be a private neighbor of $x'$ (with respect to $S_2$). 
We prove, analogously as above, that no edge 
$x'w_{x'}$ has an endvertex in $Q$, and  that  $|(H_{x'}\cup H_{w_{x'}})\cap D|\geq \gamma (H)$ 
for every $x'\in R$. 

We finish the proof of the claim by defining $$M=\{xu_x,x'w_{x'}\,|\,x\in B_1\setminus S_2, x'\in R\}\,,$$
and noting that $M$ fulfils (i) to (iv) stated above. 

Let $M'$ be the set of vertices incident to an edge of $M$, and note that one half of vertices of $M'$ are contained in $S_1\setminus S_2$ and the other half are contained in 
$S_2\setminus S_1$. If $x\in S_1\setminus (S_2\cup Q\cup M')$ and $x'\in S_2\setminus (S_1\cup Q\cup M')$, then 
every  vertex in $(\{x\}\times T_1)\cup(\{x'\}\times T_2)$ is vertically dominated (because $x\notin B_1$ and $x'\notin B_2$), and therefore 
$|(H_{x}\cup H_{x'})\cap D|\geq \gamma (H)$ (regardless if $x$ and $x'$ are adjacent or not). 
Since $|S_1|\leq |S_2|$ we have  $|S_1\setminus (S_2\cup Q\cup M')|\leq |S_2\setminus (S_1\cup Q\cup M')|$ and therefore 
$$|D|\geq |Q|\gamma (H)+|M|\gamma(H)+|(S_1\cap S_2)\setminus (B_1\cap B_2)|\gamma(H)+|S_1\setminus (S_2\cup Q\cup M')|\gamma(H).$$
Since $|S_1\setminus (S_2\cup Q\cup M')|=|S_1\setminus S_2|-\frac 12 (|Q|+|M'|)=|S_1\setminus S_2|-|B_1\cap B_2|-|M|$ 
and $|Q|=2|B_1\cap B_2|$, we conclude that $|D|\geq  |S_1|\gamma(H)\geq \gamma(G)\gamma(H)$.
\qed

\begin{corollary}\label{cor}
Let $G$ and $H$ be graphs and let $D$ be a dominating set of $G\Box H$. If there exist minimum dominating sets $D_1$ and $D_2$ of $G$ such that 
for every $h\in V(H)$ we have $p_G(D\cap G_h)\subseteq D_i$ for some $i\in [2]$, then $|D|\geq \gamma(G)\gamma(H)$.
\end{corollary}

An obvious question related to results in this paper is the question how    Theorem \ref{main} generalizes to the case 
when there are more than two pairwise distinct minimal dominating sets containing $p_G(G_h\cap D)$, when $h$ goes through $V(H)$. 
A step forward, from results obtained in the present paper, is the following conjecture (note that it is a special case of Vizing's conjecture).

\begin{conjecture}
Let $G$ and $H$ be graphs and let $D$ be a dominating set of $G\Box H$. If there exist dominating sets $S_1,S_2$ and $S_3$ of $G$ such that  for every $h\in V(H)$  we have $S_i\in M_G(p_G(D\cap G_h))$ for some  $i\in [3]$, then $|D|\geq \gamma(G)\gamma(H)$.
\end{conjecture}

 \section{Acknowledgements}

This work is supported by ARRS, grant number P1-0297.

\end{document}